\documentclass[12pt]{article}
\usepackage{amsfonts,amssymb,amsmath,theorem}
\newtheorem{Pa}{Paper}[section]
\newtheorem{Tm}[Pa]{{\bf Theorem}}
\newtheorem{La}[Pa]{{\bf Lemma}}
\newtheorem{Cy}[Pa]{{\bf Corollary}}
\newtheorem{Pn}[Pa]{{\bf Proposition}}
{
\theorembodyfont{\normalfont}
\newtheorem{Rk}[Pa]{{\bf Remark}}

}

\newcommand\norm[1]{\left\|#1\right\|}           
\newcommand\qed{\ifhmode\unskip\nobreak\fi\quad  
   \ifmmode\square\else\hbox{$\square$}\fi}      

\let\cal=\mathcal

\newcommand\xrlarrow[3][]{\genfrac{}{}{0pt}{3}%
{\displaystyle\smash[b]{\xrightarrow[#1\phantom{#3}#1]{#2}}}%
{\displaystyle\smash[t]{\xleftarrow[#1{#3}#1]{\phantom{#2}}}}}%

\numberwithin{equation}{section}

\begin{document}
\begin{center}
{\bfseries \large When a  $C^*$-algebra is a coefficient algebra\\[3pt]
for a given endomorphism}

\bigskip
V.I. Bakhtin, A.V. Lebedev

\bigskip
Belarus State University;\\
Belarus State University / University of Bialystok
\end{center}

\vspace{5mm}
\quad\parbox{14.5cm}{\small \hspace{0.5cm}
The paper presents a criterion for a $C^*$-algebra to be a
coefficient algebra\\ associated with a given endomorphism.

\medbreak
{\bfseries Keywords:} {\itshape $C^*$-algebra, endomorphism, partial
isometry, coefficient algebra,\\ transfer operator}

\medbreak
 }

\vspace{5mm}
\tableofcontents

\section{Introduction.\hfil\null\ \hbox{Coefficient algebras and transfer operators}\hfill}
\label{coefficient}

The notion of a coefficient algebra was introduced in
\cite{Leb-Odz} in connection with the study of extensions of
$C^*$-algebras by partial isometries. In this paper the authors
investigated the $C^*$-algebra $C^*({\cal A},V)$ generated by a
$C^*$-algebra  ${\cal A} \subset L(H)$ and a partial isometry $V\in L(H)$
such that the mapping ${\cal A}\ni a \mapsto VaV^*$ is an endomorphism of
${\cal A}$. It was uncovered that this investigation can be carried out
successfully if ${\cal A}$ and   $V$ satisfy the following three
conditions
\begin{gather}
\label{c1}
 Va= VaV^*V,\qquad a \in {\cal A};\\[6pt]
\label{c2}
 VaV^* \in {\cal A},\qquad a \in {\cal A}\\
\intertext{and}
\label{c3}
 V^* aV\in {\cal A},\qquad a\in {\cal A}.
\end{gather}
In \cite{Leb-Odz} the algebras  possessing these properties were
called the {\em coefficient algebras} (for $C^*({\cal A},V)$). This term
was justified by the observation that any element in $C^*({\cal A},V)$
can be presented as a Fourier like series with coefficients from
${\cal A}$ (\cite{Leb-Odz}, Theorems 2.7, 2.13). It was also shown
(\cite{Leb-Odz}, Proposition 2.2) that if ${\cal A}$ contains the
identity of $L(H)$ then the coefficient algebra can be
equivalently  defined by   conditions (\ref{c2}), (\ref{c3}) along
with the condition
\begin{equation}
\label{c4}
 V^*V\in Z({\cal A})
\end{equation}
where $Z({\cal A})$ is the center of ${\cal A}$ (it is worth mentioning that if
${\cal A}$ contains the identity then conditions (\ref{c1}), (\ref{c2})
and  (\ref{c3}) automatically imply that $V$ is a partial isometry
and the mapping ${\cal A}\ni a \mapsto VaV^*$  is an endomorphism
(\cite{Leb-Odz}, Proposition 2.2)).

In \cite{Leb-Odz} there was also presented  a certain procedure of
constructing the coefficient algebras starting from certain
initial algebras that are not coefficient ones. Following this
construction in \cite{Kw-Leb} the maximal ideal space of the
arising commutative coefficient algebras (in the situation when
the initial algebra is commutative) was described. The development
of this and other constructions, ideas and methods  led to the
general construction of covariance (crossed product type)
$C^*$-algebra for partial dynamical system \cite{Kw1}.

The interrelation between the coefficient algebras and various
crossed product type structures is not incidental.
This is due to the fact that these algebras can also be described by means
of the so-called transfer operators. Here is the description.
Using  (\ref{c4}) and recalling that  $V$ is a partial isometry
one obtains for any $a,b\in {\cal A}$ the following relations
\begin{equation}
\label{c5}
 V^*VaV^* bV = aV^*VV^*bV= aV^*bV.
\end{equation}
If we introduce the mappings
\begin{equation}
\label{c6}
 \delta \!: {\cal A}\to {\cal A}, \quad  \delta (a) =VaV^*
\end{equation}
and
\begin{equation}
\label{c7}
 \delta_* \!: {\cal A}\to {\cal A}, \quad  \delta_* (a) =V^*aV
\end{equation}
then as it has been already observed $\delta$ is an endomorphism
of ${\cal A}$ and (\ref{c5}) means that $\delta_*$ is a linear positive
mapping satisfying the condition
\begin{equation}
\label{c7'}
 \delta_*(\delta (a)b)= a\delta_* (b), \qquad a,b \in {\cal A}.
\end{equation}
The mappings $\delta_*$ possessing the latter property were called
by R. Exel \cite{exel} the {\em transfer operators} (for a given
endomorphism $\delta$). It was shown in \cite{exel} that these
objects play a significant role (they belong to the key
constructive elements) in the theory of crossed products of
$C^*$-algebras by endomorphisms.

Note that if any operator $V$ satisfies the condition
\begin{equation}
\label{c8}
 V^*VaV^* bV = aV^*bV, \qquad a,b \in {\cal A}
\end{equation}
then by passage to the adjoint we also obtain
\begin{equation}
\label{c9}
 V^*bVaV^* V = V^*bVa, \qquad a,b \in {\cal A}.
\end{equation}
Relations (\ref{c8}) and (\ref{c9}) imply in particular that
\begin{equation}
\label{c10}
aV^*V=V^*VaV^* V=V^*Va
\end{equation}
that is (\ref{c4}) is true. Thus the coefficient algebra can also
be defined as an algebra satisfying conditions (\ref{c2}),
(\ref{c3}) and (\ref{c8}). In other words a $C^*$-algebra
${\cal A}\subset L(H)$ containing the identity of $L(H)$ is the
coefficient algebra for $C^*({\cal A},V)$ iff the mapping $\delta$
(\ref{c6}) is an endomorphism and the mapping $\delta_*$
(\ref{c7}) is a transfer operator for $\delta$.

The foregoing reasoning shows on the one hand the importance of
coefficient algebras in a number of fields and on the other hand one
arrives at the natural problem: when a given $C^*$-algebra is a
coefficient algebra? The precise formulation of the problem is the
following.

{\em Let $\cal A$ be an (abstract) $C^*$-algebra containing an
identity and $\delta$ be an (abstract) endomorphism of $\cal A$.
Does there exist a triple   $(H,\pi ,U)$ consisting of a Hilbert
space~$H$, faithful non-degenerate representation $\pi\!:{\cal
A}\to  L(H)$ and a linear continuous operator  $U\!:H\to H$ such
that for every $a\in\cal A$  the following conditions are satisfied
\begin{equation}
\label{b,,4} \pi(\delta(a)) =U\pi(a)U^*,\qquad
U^*\pi(a)U \in \pi(\cal A)
\end{equation}
and
\begin{equation}
\label{c11}
U\pi(a) =\pi(\delta(a))U?
\end{equation}
That is $\pi ({\cal A})$ is the coefficient algebra for $C^*(\pi
({\cal A}), U)$ under the  fixed  endomorphism $U\,\cdot\, U^*$.}

The answer to this problem is the theme of the present article.

Since $\delta$ is an endomorphism it follows that $\delta (1)$ is
a projection and so (\ref{b,,4}) implies that $UU^*$ is a
projection, so $U$ is a partial isometry. The above discussion
means also that instead of condition  (\ref{c11}) one can use
equivalently  condition (\ref{c4}) or (\ref{c8}) (for $\pi (a)$
and $U$ respectively).

\medbreak
The article is organized as follows. In Section \ref{transfer} we
study some properties of transfer operators and in particular
we establish the conditions of existence of the so-called complete
transfer operators. By means of these operators in Section \ref{endom}
we give the answer to the problem stated above.

\section{Non-degenerate and complete transfer  operators}
\label{transfer}

Let  $\cal A$ be a $C^*$-algebra with an identity $1$ and
$\delta\!:\cal A\to \cal A$ be an endomorphism of this algebra. We
start with   recalling  some definitions and facts concerning
transfer operators.
{\mbox Every}\-thing henceforth up to Proposition \ref{Ex4.1} is
borrowed from \cite{exel}.

\medbreak
A linear map $\delta_*\!:\cal A\to \cal A$ is called a {\em transfer
operator\/} for the pair $(\cal A,\delta)$ if it is continuous and
positive an such that
\begin{equation}
\delta_*(\delta(a)b) =a\delta_*(b),\qquad a,b\in\cal A.
\label{b,,2}
\end{equation}
Since $\delta_*$ is positive it is self-adjoint an hence we have
the symmetrized version of (\ref{b,,2})
\begin{equation}
\label{b,,5}
\delta_*(b\delta(a)) =\delta_*(b)a.
\end{equation}
In particular (\ref{b,,2}) and (\ref{b,,5}) imply that $\delta_*(\cal A)$ is a two-sided ideal.

By (\ref{b,,2}) and (\ref{b,,5}) we also have that
\begin{equation}
\label{Ee2.2}
a\delta_*(1) = \delta_*(\delta(a)1) = \delta_*(1\delta(a))= \delta_*(1)a
\end{equation}
for all $a \in \cal A$, so $\delta_* (1)$ is a positive central
element in $\cal A$ and $\delta_*(1){\cal A}$ is a two-sided
ideal.

It is also worth mentioning that (\ref{b,,2}) and (\ref{b,,5}) imply
\begin{equation}
\label{Ee2.2.}
\delta_*(\,\cdot\,) = \delta_*(\delta(1)\,\cdot\,) =
\delta_*(\,\cdot\,\delta(1)).
\end{equation}

\begin{Pn}
\label{Ex2.3}
{\upshape (\cite{exel}, Proposition 2.3)}\ \
Let $\delta_*$ be a transfer operator for the pair $(\cal A,\delta)$.
Then the following are equivalent:

\smallskip
\quad\ \llap{$(i)$}\ \ the composition $E = \delta \circ \delta_*$
is a conditional expectation onto $\delta (\cal A)$,

\smallskip
\quad\ \llap{$(ii)$}\ \ $\delta\circ\delta_*\circ\delta =\delta$,

\smallskip
\quad\ \llap{$(iii)$}\ \ $\delta (\delta_*(1)) = \delta (1)$.
\end{Pn}

If the equivalent conditions of Proposition  \ref{Ex2.3} hold then
R.\ Exel calls $\delta_*$ a {\em non-degenerate} transfer operator.

\begin{Pn}\label{Ex2.5}
{\upshape (\cite{exel}, Proposition 2.5)}\ \
Let $\delta_*$ be a non-degenerate transfer operator. Then $\cal A$
may be written as the direct sum of ideals
\[
{\cal A} = \mathop{\mathrm{Ker}}\delta \oplus
\mathop{\mathrm{Im}} \delta_*
\]
\end{Pn}

\begin{Pn} \label{Ex4.1}
{\upshape (\cite{exel}, Proposition 4.1)}\ \
$\delta ({\cal  A})$ is a hereditary subalgebra of $\cal A$ iff \
$\delta ({\cal A})= \delta (1){\cal A}\delta (1)$.
\end{Pn}

\begin{Rk}\label{nonunique}
If $\delta $ is given then in general a non-degenerate transfer
operator $\delta_*$ (if it exists) is not defined in a unique way.

\medbreak
{\em Example.} Let  ${\cal A} = C(X)$ where $X = {\mathbf R}\,
(\mathrm{mod}\, 1)$ and let $\delta \!: {\cal A} \to {\cal A}$ be given by
the formula
\[
\delta (a) (x) = a (2x\, (\mathrm{mod}\, 1)).
\]
Take any continuous function $\rho$ on $X$ having the properties
\begin{gather*}
0\le \rho (x) \le 1, \qquad x\in X,\\[6pt]
\rho \left(x+\frac{1}{2}\right) = 1 -\rho (x), \qquad
0\le x \le \frac{1}{2}.
\end{gather*}
Let us define the operator $\delta_*$ on $C(X)$ by the formula
\[
\delta_*  (a)(x) =   a\left(\frac{x}{2}\right) \rho
\left(\frac{x}{2}\right) + a\left(\left[\frac{x}{2} +
\frac{1}{2}\right]\right)\rho \left(\left[\frac{x}{2} +
\frac{1}{2}\right]\right),
\]
where $[x] = x \, (\mathrm{mod}\, 1)$.
Clearly for any  $\rho $ chosen $\delta_*$ is a transfer operator
for $\delta$ and since $\delta_*(1)=1$ it is non-degenerate.
\end{Rk}

Now we note that the non-degeneracy of transfer operators also implies
some additional properties.
\begin{Pn}\label{l1}
Let $\delta_*$ be a non-degenerate transfer operator. Then

\smallskip
\quad\llap{$1)$}\ \ $\delta_* (1)$ is an orthogonal projection lying
in the center of $\cal A$,

\smallskip
\quad\llap{$2)$}\ \ $\delta_* ({\cal A}) = \delta_* (1){\cal A}$.

\smallskip
\quad\llap{$3)$}\ \ $\delta_* \!: \delta ({\cal A}) \to \delta_*
({\cal A})$ is a $^*$-isomorphism $($the inverse to the $^*$-isomorphism
$\delta \!: \delta_*({\cal A}) \to \delta({\cal A})$$)$.
\end{Pn}

{\bfseries Proof.} 1) \, It has been already observed that $\delta_*
(1)$ is a central element in $\cal A$ so it is enough to show that
it is a projection.

By Proposition  \ref{Ex2.5} we have that
\[
{\cal A} = \mathop{\mathrm{Ker}}\delta \oplus \mathop{\mathrm{Im}}
\delta_*.
\]
Thus the mapping
\begin{equation}\label{*1}
\delta \!: \mathop{\mathrm{Im}} \delta_* \to \delta({\cal A})
\end{equation}
is a $^*$-isomorphism and so the mapping
\begin{equation}\label{***}
\delta^{-1} \!: \delta({\cal A}) \to \mathop{\mathrm{Im}} \delta_*
\end{equation}
is a $^*$-isomorphism as well.

Since $\delta_*$ is non-degenerate it follows from (iii) of
Proposition  \ref{Ex2.3}  that
\begin{equation}
\label{**1}
\delta (\delta_*(1)) = \delta (1).
\end{equation}
Now (\ref{***}) and (\ref{**1}) imply
\[
\delta_* (1) = \delta^{-1}(\delta (1)).
\]
Since $\delta (1)$ is a projection (as $\delta$ is an endomorphism)
and $\delta^{-1}$ is a morphism it follows that $\delta^{-1}(\delta
(1))$ is a projection as well. So 1) is proved.

\smallbreak
2) \, Using (\ref{b,,5}) we obtain
\begin{equation}\label{***.}
\delta_*(1){\cal A} = \delta_*(1\delta({\cal A}))
\subset  \delta_*({\cal A}).
\end{equation}
In view of the non-degeneracy of $\delta_*$ and  (iii) of
Proposition  \ref{Ex2.3} we have
\begin{equation}\label{***..}
\delta (\delta_*(1){\cal A}) = \delta
(\delta_*(1))\cdot \delta({\cal A}) = \delta (1) \delta({\cal A})
= \delta({\cal A}).
\end{equation}
Now we conclude from (\ref{***.}), (\ref{*1}) and (\ref{***..}) that
the mapping
\begin{equation}
\label{***...} \delta \!: \delta_* (1) {\cal A} \to \delta
({\cal A})
\end{equation}
is a $^*$-isomorphism (we have already noticed that $\delta_* (1)
{\cal A}$ is an ideal (recall (\ref{Ee2.2}))). Therefore
(\ref{***...}) and (\ref{*1}) imply the equality $\delta_*
(1){\cal A}= \delta_* ({\cal A})$ and the proof of 2)  is
finished.

\smallbreak
3) \,For any $a\in \cal A$ we have that $\delta(a)\in\delta(\cal A)$ and $\delta_*\delta(a)\in
\delta_*(\cal A)$, and $\delta\delta_*\delta(a) =\delta(a)$ by Proposition \ref{Ex2.3} (ii).
Therefore $\delta_*\!:\delta(\cal A)\to \delta_*(\cal A)$ is a right inverse to $\delta\!:
\delta_*(\cal A)\to \delta(\cal A)$. But since it has been already observed that
$\delta\!: \delta_*(\cal A)\to \delta(\cal A)$ is an isomorphism it follows that $\delta_*$ is an
isomorphism as well.
\qed

\begin{Pn}\label{almost-unique}
Let $\cal A$ be a $C^*$-algebra with an identity $1$,
\,$\delta \!: {\cal A} \to {\cal A}$ be an endomorphism of $\cal A$
and $\delta_{*i}, \, i=1,2$ be two non-degenerate transfer operators
for $({\cal A}, \delta)$. Then

\smallskip
\quad\llap{$1)$}\quad $\delta_{*1} (1)= \delta_{*2} (1)$,

\smallskip
\quad\llap{$2)$}\quad $\delta_{*1}({\cal A}) = \delta_{*2}({\cal A})$.
\end{Pn}

{\bfseries Proof.} \, By  1) of Proposition  \ref{l1} \, we have that
$\delta_{*i}(1), \, i=1,2$ are the orthogonal projections
belonging to the center of $\cal A$. Set
\[
P = \delta_{*1} (1)\delta_{*2} (1).
\]
By 1) of Proposition \ref{l1} we have
\begin{equation}
\label{i}
P \in \delta_{*i} ({\cal A}),\quad \ i=1,2,
\end{equation}
and by the non-degeneracy of $\delta_{*i}$, $i=1,\,2$ (see (iii) of
Proposition \ref{Ex2.3}) we obtain
\begin{equation}\label{ii}
\delta (P)= \delta (\delta_{*1} (1)\delta_{*2}(1))= \delta
(\delta_{*1} (1))\delta (\delta_{*2} (1))= \delta (1)\delta (1)
=\delta (1).
\end{equation}
But in view of  Proposition  \ref{Ex2.5}  the mappings
\begin{equation}\label{*1*}
\delta \!:  \delta_{*i}({\cal A}) \to \delta({\cal A}), \quad i=1,2
\end{equation}
are $^*$-isomorphisms.
Now (\ref{i}), (\ref{ii}) and (\ref{*1*}) imply
\[
P = \delta_{*1} (1)=\delta_{*2} (1).
\]

Since by 1) $\delta_{*1} (1)=\delta_{*2} (1)$ we conclude by 2)
of Proposition \ref{l1} that
 $\delta_{*1} ({\cal A}) = \delta_{*2}({\cal A})$.
The proof is complete.
\qed

\begin{Rk}
It has been observed that in general a non-degenerate transfer
operator (for a given pair $({\cal A}, \delta)$)
is not unique. Part 2) of Proposition  \ref{almost-unique} and part 3) of
Proposition  \ref{l1}
tell us that nevertheless
its restriction onto    $\delta ({\cal A})$ is unique.
\end{Rk}

We shall call the transfer operator $\delta_*$  {\em complete,} if
\begin{equation}
\label{b,,3}
\delta\delta_*(a) =\delta(1)a\delta(1),\qquad a\in\cal A.
\end{equation}
Observe that a complete transfer operator is non-degenerate. Indeed,
(\ref{b,,3}) implies
\[
\delta\delta_*\delta(a) = \delta(1)\delta (a) \delta(1) = \delta (a)
\]
and so condition (ii) of Proposition  \ref{Ex2.3} is satisfied.

Note also that (\ref{b,,3})  implies
\begin{equation}
\label{d*dd*}
\delta_{*}\delta\delta_* =\delta_* .
\end{equation}
The next result presents the criteria for the existence of a complete
transfer operator.

\begin{Tm}\label{complete}
Let $\cal A$ be a $C^*$-algebra with an identity $1$ and $\delta \!:
{\cal A} \to {\cal A}$ be an endomorphism of $\cal A$. The following
are equivalent:

\medskip
\quad\llap{$1)$}\ \
there exists a complete transfer operator $\delta_*$ $($for
$({\cal A}, \delta )$$),$

\medskip
\quad\llap{$2)$}\qquad\llap{$(i)$}\ \,
there exists a non-degenerate transfer operator $\delta_*$ and

\smallskip
\quad\qquad\llap{$(ii)$}\ \,
$\delta ({\cal A})$ is a hereditary subalgebra of $\cal A$$;$

\medskip
\quad\llap{$3)$}\qquad\llap{$(i)$}\ \,
there exists a central orthogonal projection  $P\in \cal A$ such that

\smallskip
\qquad\qquad\quad\llap{$a)$}\ $\delta (P) = \delta (1)$,

\smallskip
\qquad\qquad\quad\llap{$b)$}\ the mapping $\delta \!: P {\cal A} \to
\delta ({\cal A})$ is a $^*$-isomorphism, and

\smallskip
\quad\qquad\llap{$(ii)$}\ \, $\delta ({\cal A}) = \delta (1){\cal A}
\delta (1)$.

\medbreak\noindent
Moreover the objects in $1)$ -- $3)$ are defined in a unique way
$($i.\,e.\ the transfer operator $\delta_*$ in $1)$ and $2)$ is
unique and the projection $P$ in $3)$ is unique as well$)$ and
\begin{equation}
\label{P} P=\delta_*(1)
\end{equation}
and
\begin{equation}\label{d*-}
\delta_*(a) =\delta^{-1}(\delta(1)a\delta(1)), \ \  a\in {\cal A}
\end{equation}
where $\delta^{-1}\!:\delta({\cal A})\to P\cal A$ is the inverse
mapping to\ \ $\delta\!:P\cal A\to \delta(\cal A)$.
\end{Tm}

{\bfseries Proof.} \, 1) $\Rightarrow$ 2). \, By Proposition
\ref{Ex4.1} \ \ 2) (ii) and 3) (ii) are equivalent.

We have already noticed that if $\delta_*$ is complete then it is
non-degenerate. Moreover (\ref{b,,3}) implies
\[
\delta ({\cal A}) \supset \delta (1){\cal A} \delta (1)
\]
and on the other hand
\[
\delta ({\cal A}) = \delta (1)\delta ({\cal A}) \delta (1)\subset
\delta (1){\cal A} \delta (1).
\]
Thus 1) implies 2).

\medbreak
2) $\Rightarrow$ 3). \, We have already noticed that 2) (ii) and 3)
(ii) are equivalent.

Set $P = \delta_* (1)$. By 1) of Proposition \ref{l1} $P$ is an
orthogonal central projection, and by (iii) of Proposition \ref{Ex2.3}
$\delta (P) = \delta (1)$. So 3) (i) (a) is true.

By 2) of Proposition  \ref{l1} we have that $\delta_* ({\cal A}) =
P\cal A$ and it follows from Proposition \ref{Ex2.5} that the
mapping $\delta \!: P {\cal A} \to \delta ({\cal A})$ is a
$^*$-isomorphism.
So 3) (i) (b) is true.

\medbreak
3) $\Rightarrow$ 1). \, Let $\delta^{-1}\!:\delta({\cal A})\to
P\cal A$ be the inverse mapping to\ \ $\delta\!:P\cal A\to
\delta(\cal A)$. Define the operator  $\delta_*$ by the formula
$\delta_*(a) =\delta^{-1}(\delta(1)a\delta(1))$. Clearly
$\delta_*$ is positive and satisfies (\ref{b,,3}). Note that
\begin{equation*}
\delta\bigl(\delta_*(\delta(a)b)\bigr) =\delta(1)\delta(a)b
\delta(1) =\delta(a)\delta(1)b\delta(1) =\delta(a\delta_*(b)).
\end{equation*}
But in view of the definition of  $\delta_*$ the elements  $\delta_*(\delta(a)b)$
and
$a\delta_*(b)$ belong to the ideal  $P\cal A$ where the endomorphism $\delta$
is injective. Therefore they coincide and thus  (\ref{b,,2}) is proved. So  $\delta_*$
is a complete transfer operator.

\medbreak
Now to finish the proof let us verify the uniqueness of the objects
mentioned in 1)~--~3).

The uniqueness of the projection $P$ in 3) is in fact established
in 1) of Proposition~\ref{almost-unique} (since here we have $P
= \delta_* (1)$ for some non-degenerate $\delta_*$ which also
proves~(\ref{P})).

Recalling (\ref{Ee2.2.}) we obtain
\begin{equation}\label{+}
\delta_*(a) = \delta_*(\delta(1)a\delta(1)), \ \, a\in \cal A .
\end{equation}
But by 3) (ii)\ \ $\delta(1)a\delta(1) \in \delta ({\cal A})$.
Therefore (\ref{+}) and 3) of Proposition \ref{l1} imply the
uniqueness of $\delta_*$.

Finally we note that formula   (\ref{d*-}) has been
already proven in the course of the proof of 3) $\Rightarrow$ 1).

The proof is complete.
\qed

\medbreak
In connection with Theorem \ref{complete} it makes sense to note the
next useful observation which is given in Proposition \ref{partial
aut}. We recall that a {\em partial automorphism} of a
\hbox{$C^*$-alge}\-bra $\cal A$ is a triple $(\theta , I,J)$ where $I$
and $J$ are closed two-sided ideals of $\cal A$ and $\theta \!: I \to
J$ is a $^*$-isomorphism.

\begin{Pn}\label{partial aut}
Let $\delta_*$ be a complete transfer operator
for $({\cal A}, \delta )$ and $\delta (1)\in Z({\cal A})$, then

$1)$ \, both the triples \, $(\delta , \delta_* (1){\cal A}, \delta
(1){\cal A})$ \, and \, $(\delta_* , \delta (1){\cal A}, \delta_*
(1){\cal A})$ \, are partial automorphisms (that are inverses to
each other), and

$2)$ \, $\delta_* \!: {\cal A} \to {\cal A}$ is an endomorphism.
\end{Pn}

{\bfseries Proof.}\, 1) By Proposition \ref{l1} we have that $\delta_*
(1)$ is a central projection. Since by the condition of the lemma
$\delta (1)$ is a central projection as well it follows that both
$\delta_* (1){\cal A}$ and $\delta (1){\cal A}$ are the ideals. In
addition by 3) (ii) of Theorem \ref{complete} we have that $\delta
({\cal A})=\delta (1){\cal A}\delta (1) =\delta (1){\cal A}$. Now
1) follows from Proposition \ref{l1}.

\smallbreak
2) Using (\ref{d*-}) and the condition that $\delta (1)\in
Z({\cal A})$ and recalling that $\delta_* \!: \delta ({\cal A}) \to
\delta_*({\cal A})$ is a $^*$-isomorphism we obtain for any $a,b
\in \cal A$
\begin{gather*}
\delta_*(ab)=\delta^{-1}(\delta(1)ab\delta(1))=
\delta^{-1}(\delta(1)a\delta(1)\delta(1)b\delta(1))\\[6pt]
=\delta^{-1}(\delta(1)a\delta(1)) \,
\delta^{-1}(\delta(1)b\delta(1))= \delta_*(a)\delta_*(b).
\end{gather*}
Thus $\delta_* \!: {\cal A} \to {\cal A}$ is an endomorphism.
\qed

\section{Coefficient algebras associated with  given\hfil\null\ \hbox{endomorphisms}\hfill}
\label{endom}

Now we pass to the main result  of the article.

Let $\cal A$  be a $C^*$-algebra  containing an identity $1$ and
$\delta\!:\cal A\to\cal A$ be an endomorphism. We say that $\cal
A$ is a {\em coefficient algebra associated with $\delta$} if
there exists a triple $(H,\pi ,U)$ consisting of a Hilbert space
$H$, faithful non-degenerate representation $\pi\!:{\cal A}\to
L(H)$ and a linear continuous operator $U\!:H\to H$ such that
conditions (\ref{b,,4}) and (\ref{c11}) are
satisfied.

\begin{Tm}\label{b..1}
$\cal A$ is a  coefficient algebra associated with
$\delta$ iff there exists a complete transfer operator $\delta_*$
for $({\cal A}, \delta )$ (that is either of the equivalent
conditions of Theorem \ref{complete} hold).
\end{Tm}

{\bfseries Proof.} \, {\em Necessity}. \, If conditions (\ref{b,,4}) and
(\ref{c11})  are satisfied then (identifying $\cal A$ with
$\pi (\cal A)$) one can set
\[
\delta_* (\cdot ) = U^* (\cdot )U
\]
and it is easy to verify that $\delta_*$ is a complete transfer
operator.

\smallbreak
{\em Sufficiency.} \, Let $\delta_*$ be a complete transfer operator.
We shall construct the desired Hilbert space $H$ by means of the
lements of the initial algebra $\cal A$ in the following way. Let
$\langle\,\cdot\,,\, \cdot\,\rangle$ be a certain non-negative inner
product on $\cal A$ (differing from a common inner product only in
such a way that for certain non-zero elements $v\in\cal A$ the
expression $\langle v,v\rangle$ may be equal to zero). For example
this inner product may have the form $\langle v,u\rangle =f(u^*v)$
where $f$ is some positive linear functional on $\cal A$. If one
factorizes $\cal A$ by all the elements $v$ such that $\langle
v,v\rangle =0$ then he obtains a linear space with a strictly positive
inner product. We shall call the completion of this space with respect
to the norm $\norm v =\sqrt{\langle v,v \rangle}$ the {\em Hilbert
space generated by the inner product\/} $\langle\,\cdot\,,\,
\cdot\,\rangle$.

Let $V$ be the set of all positive linear functionals on $\cal A$. The
space $H$ will be the completion of the direct sum $\bigoplus_{f\in
V}H^f$ of some Hilbert spaces $H^f$. Every $H^f$ will in turn be the
completion of the direct sum of Hilbert spaces $\bigoplus_{n\in\mathbb
Z} H^f_n$. The spaces $H^f_n$ are generated by non-negative inner
products $\langle\,\cdot\,,\,\cdot\,\rangle_n$ on the initial algebra
$\cal A$ that are given by the following formulae
\begin{align}
\label{b,,7}
\langle v,u\rangle_0 &=f(u^*v);\\[3pt]
\label{b,,8}
\langle v,u\rangle_n &=f\bigl(\delta_*^n(u^*v)\bigr),
\qquad n\ge 0;\\[3pt]
\label{b,,9}
\langle v,u\rangle_n
&=f\bigl(u^*\delta^{|n|}(1)v\bigr), \qquad n\le 0.
\end{align}
The properties of these inner products are described in  the next

\begin{La}
\label{b..2}
For any \/ $v,u\in\cal A$ the following equalities are true
\begin{align}
\label{b,,10} \langle \delta(v),u\rangle_{n+1} &=\langle
v,\delta_*(u)\rangle_n, \qquad n\ge 0;\\[3pt]
\label{b,,11} \langle\delta^{|n|}(1)v,u\rangle_{n+1} &=\langle
v,\delta^{|n|}(1)u \rangle_n,\qquad n<0.
\end{align}
\end{La}

{\bfseries Proof.}
Indeed, the proof of  (\ref{b,,10}) reduces to the verification of the equality
\[
\delta_*^{n+1}(u^*\delta(v)) =\delta_*^n (\delta_*(u^*)v)
\]
which follows from  (\ref{b,,5}), and the proof of
(\ref{b,,11}) reduces to the verification of the equality
\[
u^*\delta^{|n|-1}(1)\delta^{|n|}(1)v = u^*\delta^{|n|}(1)
\delta^{|n|}(1)v
\]
which follows from the equalities
\[
\delta^{|n|-1}(1) \delta^{|n|}(1) =\delta^{|n|-1}(1\cdot\delta(1))
=\delta^{|n|}(1) =(\delta^{|n|}(1))^2. \qed
\]

\smallbreak
Now let us define the operators  $U$ and $U^*$ on the space  $H$
constructed. These operators will leave invariant all the
subspaces  $H^f\subset H$. The action of  $U$ and $U^*$ on every
$H^f$ is the same and its scheme is presented in the first line
of the next diagram.

\noindent
\begin{alignat*}{14}
&\dots\ \,&&
\xrlarrow[\,]{\delta^3(1)\,\cdot\,}{\delta^3(1)\,\cdot\,}&&
\ H^f_{-2}\ &&
\xrlarrow[\,]{\delta^2(1)\,\cdot\,}{\delta^2(1)\,\cdot\,}&&
\ H^f_{-1}\ &&
\xrlarrow[\,]{\delta(1)\,\cdot\,}{\delta(1)\,\cdot\,}&&
\ H^f_{0}\ &&
\xrlarrow[\,]{\delta}{\delta_*}&&
\ H^f_{1}\ &&
\xrlarrow[\,]{\delta}{\delta_*}&&
\ H^f_{2}\ &&
\xrlarrow[\,]{\delta}{\delta_*}&
\ \,\dots&&&
\qquad \xrlarrow[\ ]{\textstyle U}{\textstyle\,U^*}\\[6pt]
&\dots &&&&\,\mspace{0mu}\delta^2(a) &&&&\ \delta(a) &&&&\ \;a
&&&&\ \;a &&&&\ \;a &&&
\dots &&&\qquad\ \:\pi(a)
\end{alignat*}

\medbreak\noindent Formally this action is defined in the
following way. Consider any {\em finite} sum
\[
h =\bigoplus_n h_n \in H^f, \ \ h_n\in H^f_n .
\]
Set
\[
Uh =\bigoplus_n(Uh)_n \qquad\text{and}\qquad U^*h =\bigoplus_n(U^*h)_n
\]
where
\begin{align}
\label{b,,12}
(Uh)_n &=
\begin{cases}
\delta(h_{n-1}),&\ \ {\rm if}\ \ n>0,\\
\delta^{|n|+1}(1)h_{n-1},& \ \ {\rm if}\ \ n\le 0,
\end{cases}\\[6pt]
(U^*h)_n &=
\begin{cases}
\delta_*(h_{n+1}), & \ \ {\rm if}\ \ n\ge 0,\\
\delta^{|n|}(1)h_{n+1},&\ \  {\rm if}\ \ n<0.
\end{cases}
\label{b,,13}
\end{align}
Lemma  \ref{b..2} guarantees that the operators  $U$ and $U^*$ are
well defined (i.\,e.\ they preserve factorization and completion by
means of which the spaces $H^f_n$ were built from the algebra~$\cal A$)
and $U$ and $U^*$ are mutually adjoint.

Now let us define the representation  $\pi\!:{\cal A}\to L(H)$.
For any $a\in\cal A$ the operator  $\pi(a)\!:H\to H$ will leave
invariant all the subspaces $H^f\subset H$ and also all the
subspaces  $H^f_n \subset H^f$. If  $h_n\in H^f_n$ then we set
\begin{equation}
\label{b,,14}
\pi(a)h_n =
\begin{cases}
ah_n, & \ \ n\ge 0,\\
\delta^{|n|}(a)h_n, &\ \ n\le 0.
\end{cases}
\end{equation}
The scheme of the action of the operator  $\pi(a)$ is presented in
the second line of the diagram given above.

Let us verify  equalities  (\ref{b,,4}) for the representation $\pi$.
Take any   $h_n \in H^f_n$. Then for  $n<0$ we have
\begin{gather*}
U^*\pi(a)Uh_n =\delta^{|n|}(1)\delta^{|n|-1}(a)\delta^{|n|}(1)h_n,\\[6pt]
\pi(\delta_*(a))h_n =\delta^{|n|}(\delta_*(a))h_n =
\delta^{|n|}(1)\delta^{|n|-1}(a)\delta^{|n|}(1)h_n
\end{gather*}
(where the final equality follows from  (\ref{b,,3}));
and for  $n\ge 0$ we have
\begin{gather*}
\pi(\delta_*(a))h_n =\delta_*(a)h_n,\\[6pt]
U^*\pi(a)Uh_n =\delta_*(a\delta(h_n)) =\delta_*(a)h_n.
\end{gather*}
In addition for $n\le 0$ one has
\begin{gather*}
\pi(\delta(a))h_n =\delta^{|n|+1}(a)h_n,\\[6pt]
U\pi(a)U^*h_n =\delta^{|n|+1}(1)\delta^{|n|+1}(a)\delta^{|n|+1}(1)
h_n =\delta^{|n|+1}(a)h_n;
\end{gather*}
and for  $n>0$
\begin{gather*}
\pi(\delta(a))h_n =\delta(a)h_n,\\[6pt]
U\pi(a)U^*h_n =\delta(a\delta_*(h_n)) =\delta(a)h_n\delta(1),
\end{gather*}
and moreover (\ref{b,,8}) and (\ref{b,,5}) imply that for $n>0$ the
element $\delta(a)h_n\delta(1)$ coincides with $\delta(a)h_n$ in the
space $H^f_n$.

Thus we have proved that $U\pi(a)U^* =\pi(\delta(a))$ and
$U^*\pi(a)U =\pi(\delta_*(a))$ for any $a\in\cal A$.

To finish the proof it is enough to observe the faithfulness of the
representation $\pi$. But this follows from the definition of the
inner product in (\ref{b,,7}),
the definition of $\pi$ (see the second
line in the diagram) and the standard Gelfand-Naimark faithful
representation of a $C^*$-algebra.
The proof is complete.
\qed

\medbreak
We shall say that an  {\em endomorphism $\delta \!: \cal A\to \cal
A$ is generated by an isometry } if there exists a triple
$(H,\pi ,U)$ consisting of a Hilbert space  $H$, faithful
non-degenerate representation $\pi\!:{\cal A}\to  L(H)$ and an
{\em isometry} $U\!:H\to H$ such that for every
$a\in\cal A$  the following conditions are satisfied
\begin{equation}
\label{b''4} \pi(\delta(a)) =U\pi(a)U^*,\qquad U^*\pi(a)U \in \pi
(\cal A).
\end{equation}
Theorem \ref{b..1} implies in particular an obvious

\begin{Cy}
An endomorphism $\delta$ is generated by an  isometry iff there
exists a complete transfer operator $\delta_*$ for $({\cal A},
\delta )$  such that $\delta_*(1)=1$ or, equivalently iff $\delta$
is a monomorphism with hereditary range.
\end{Cy}

{\bf Proof.} \, Evidently $\delta$ is generated by an  isometry iff
$\cal A$ is a  coefficient algebra associated with $\delta$ and
 we have $U^*U=1$. So in this case $\delta_*(1)=1$ and
the projection $P$ mentioned in 3) of Theorem \ref{complete} is
equal to $1$.
\qed

\end{document}